\documentclass[12pt]{amsart}
\usepackage{color}
\usepackage{amssymb}
\usepackage{graphicx}

\numberwithin{equation}{section}

\renewcommand{\Im}{\mathop{\rm Im}\nolimits}
\newcommand{\RR}{\mathbb R}
\newcommand{\ZZ}{\mathbb Z}
\newcommand{\defeq}{\stackrel{\rm{def}}{=}}

\begin{document}

\thispagestyle{empty}
\title[Effective dynamics for $ N$-solitons]{Effective dynamics for $N$-solitons of the Gross-Pitaevskii equation}

\author{Trevor Potter}
\address{Department of Mathematics, University of California,
Berkeley, CA  94720, USA}
\email{potter@math.berkeley.edu}

\maketitle

\begin{abstract}
We consider several solitons moving in a slowly
varying external field. We show that the effective dynamics
obtained by restricting the full Hamiltonian to the 
finite dimensional manifold of $ N$-solitons (constructed
when no external field is present) provides a remarkably
good approximation to the actual soliton dynamics. That
is quantified as an error of size $ h^2 $ where $ h $ is
the parameter describing the slowly varying nature of the
potential. This also indicates that previous mathematical results
of Holmer-Zworski \cite{holmer-zworski} 
for one soliton are optimal. For potentials
with unstable equilibria the Ehrenrest time, $ \log(1/h)/h $, appears
to be the natural limiting time for these effective dynamics. 
We also show that the results of Holmer-Perelman-Zworski 
\cite{holmer-zworski1} for two mKdV solitons apply numerically
to a larger number of interacting solitons.
We illustrate
the results by applying the method with the external potentials 
used in Bose-Einstein soliton train experiments of 
Strecker \emph{et al} \cite{strecker}.
\end{abstract}

\section{Introduction}
\label{sec:intro}

In many situations a wave moving in a slowly varying field,
that is, a field described by a potential whose derivatives 
are much smaller than the oscillations/width of the wave, 
can be described using classical dynamics. This is the 
basis of the semiclassical/short wave approximation, perhaps best
known in the case of the linear Schr\"odinger equation,
\begin{equation}
\label{eq:SE}
  i h \partial_tu = -\frac{1}{2} h^2 \partial_x^2 u + V( x)u \,,  
\end{equation}
where $ V $ is an infinitely differentiable potential. 
A typical result concerns a propagation of a coherent state
\[   u ( x , 0 ) = \exp \left( \frac i h 
\left( ( x - x_0 ) \xi_0 + i (x-x_0)^2 /2
\right) \right) \,, \]
maximally concentrated near the point $ ( x_0 , \xi_0 ) $ in the 
position-momentum space. In that case,
\begin{equation}
\label{eq:cs} u ( x , t ) =   a_0 ( x, t ) 
\exp \left(  \frac i h \varphi ( x , t ) \right) 
 + {\mathcal O} ( h^{\frac 12} ) \,, \ \  0 < t < T ( h )   \,, 
\end{equation}
where $ \Im \partial_x^2\varphi  > 0 $, $ \Im \varphi \geq 0 $, and 
\[ \Im \varphi ( x, t ) = 0 \ \Rightarrow \ x = x ( t ) \,, \ \ 
\partial_x \varphi ( x, t ) = \xi ( t ) \,, \]
where $  ( x(t) , \xi ( t ) ) $ satisfy Newton's
equations:
\begin{equation}
\label{eq:Newt}  x' ( t ) = \xi( t) \,, \ \ \xi' ( t) = -  V' ( x ) \, , \ \ 
x( 0 ) = x_0 \,, \ \ \xi ( 0 ) = \xi_0 \,.
\end{equation}
The time of the validity of \eqref{eq:cs}, $ T ( h ) $, depends on 
the properties of the flow \eqref{eq:Newt}, and in general it is 
limited by the {\em Ehrenfest time},
\begin{equation}
\label{eq:Ehr}  T ( h )  \sim \ \log \left( \frac 1 h \right) \,,
\end{equation}
see \cite{DR} for a recent discussion on the case of one dimension.

\begin{figure}[t]
\begin{tabular}{c c}
\includegraphics[scale=.5]{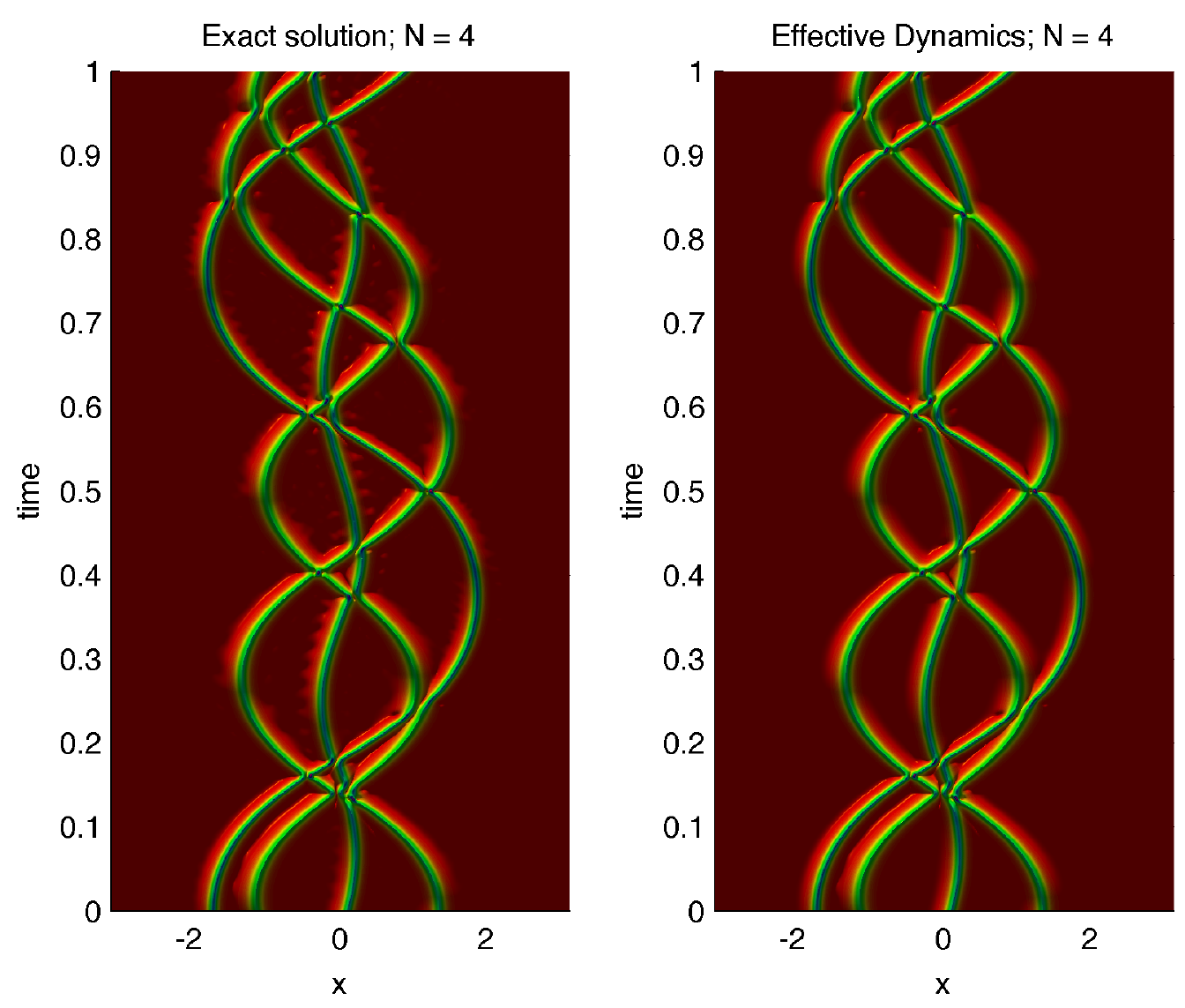} & \includegraphics[scale=.5]{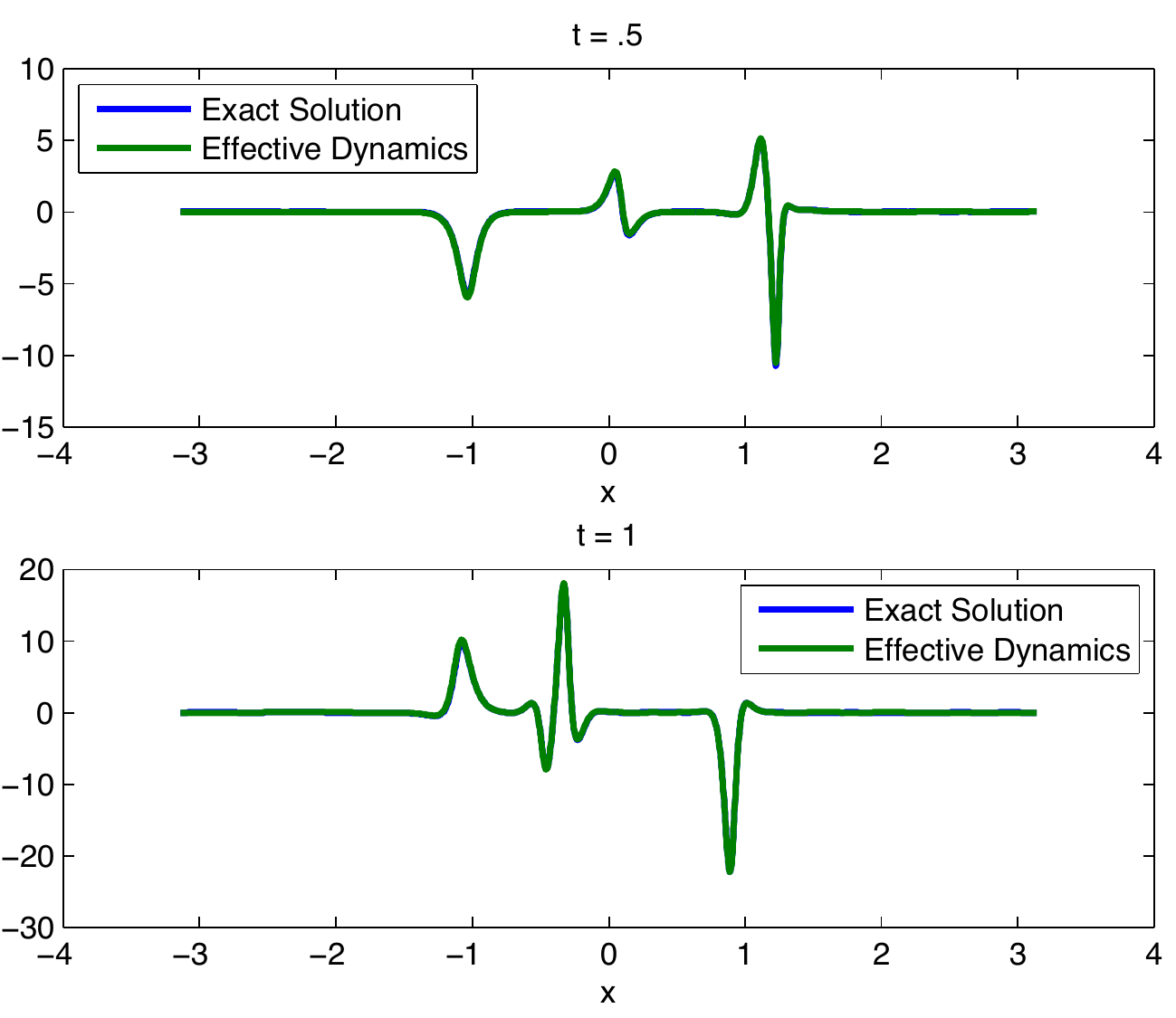} \\
\end{tabular}
\caption{A side-by-side comparison of the effective dynamics versus the exact solution
of (\ref{eq:GP}) for 4 solitons with the potential $W(x) = -100e^{\cos x}$. 
The plot on the left shows the absolute value of the solutions up to time $t = 1$.
The plot on the right shows the real part of the solutions at times $t = 0.5$ and
$t = 1$.
Compared to the solutions in Figures \ref{fig:comp} and \ref{fig:keepshape},  
much less discrepancy between the two solutions is visible.}
\label{fig:perfectmatch}
\end{figure}

The approximation (\ref{eq:cs}) means that the solution is 
concentrated for logarithmically 
long times on classical trajectories. 
The phase $ \varphi $ and the amplitude
$ a_0 $ 
can be described very precisely and $ a_0 $ can be refined to give an 
asymptotic expansion -- 
see \cite{Ha} for an early mathematical treatment 
and \cite{Ro} for more recent developments and references.

In this paper we consider the Gross-Pitaevski equation,
which is the cubic non-linear 
Schr\"odinger equation with a potential:
\begin{equation}
\label{eq:GP0}
ih \partial_tu = -\frac{1}{2} h^2 \partial_x^2 u - u|u|^2 + V(x)u\,.
\end{equation}
It provides a mean field approximation for the evolution of Bose-Einstein
condensate in an external field given by the potential $ V ( x ) $ 
-- see the monograph \cite{PS} and references given there. 
Questions about propagation of localized states are also
natural in the setting of \eqref{eq:GP0} and have been much 
studied. One direction is described in a recent monograph \cite{Ca}.

In this note we present a numerical 
study of multiple soliton propagation
for \eqref{eq:GP0} and show that it can be described 
very accurately using a natural effective dynamics 
-- see Figure \ref{fig:perfectmatch}. That effective
dynamics is based on mathematical results of 
Holmer-Zworski \cite{holmer-zworski} and Holmer-Perelman-Zworski
\cite{holmer-zworski1} and we refer to those papers for
pointers to earlier mathematical works on that subject. 

Following the convention of earlier papers -- see Fr\"ohlich et al 
\cite{FrSi} -- 
we rescale equation \eqref{eq:GP0}
so that the parameter $ h $ is in the potential which is now slowly 
varying:
\begin{equation}
\label{eq:GP}
i \partial_tu = -\frac{1}{2} \partial_x^2 u - u|u|^2 + V(x)u\,, \qquad V(x) = W(hx)\,.
\end{equation}
For $ V \equiv 0 $ this equation is completely integrable 
--  see for instance \cite{FT}. 
One of the most striking consequences of that is the existence
of exact $ N $-soliton solutions:
\begin{gather}
\label{eq:uxt0}
\begin{gathered}
   u ( x, t ) = q_N ( x , a + t v, v , \theta + \frac{t}{2}( \mu^2  + v^2 ), \mu) \,,
\\ a, v  \in \RR^N \,, \ \ \theta \in (\RR/ 2 \pi \ZZ )^N \,, \ \ 
\mu \in \RR^N_+ \,, 
\end{gathered}
\end{gather}
where the construction of $ q_N = q_N ( x , a , v, \theta, \mu ) $ 
will be recalled in \S \ref{sec:Ns}.

When $ V \not \equiv 0 $ and 
\[ u (x , 0 ) = q_N ( x , a, v, \theta, \mu ) \,, \]
the exact dynamics \eqref{eq:uxt0} 
is replaced by 
\begin{gather}
\label{eq:uxtV}
\begin{gathered}
   u ( x, t ) = q_N ( x, a(t), v ( t) , \theta ( t ) , \mu (t)) + 
{\mathcal O} ( h^2 ) \,, 
\end{gathered}
\end{gather}
where the precise meaning of the error, and its optimality,
will be described below.
The parameters of the multisoliton approximation solve the system
of ordinary differential equations:
\begin{gather}
\label{eq:effd}
\begin{gathered}
\dot{v_j} = -\mu_j^{-1} (\partial_{a_j} V_N + v_j \partial_{\theta_j} V_N)\,, \ \ 
\dot{a}_j  = v_j + \mu_j^{-1} \partial_{v_j}V_N\,,\\
\dot{\mu}_j  = \partial_{\theta_j} V_N \,, \ \ 
\dot{\theta}_j  = v_j^2/2 + \mu^2_j/2 + \mu_j^{-1} v_j \partial_{v_j} V_N - \partial_{\mu_j} V_N \,,
\end{gathered}
\end{gather}
and where
\[
 V_N(a,v,\theta,\mu) \stackrel{\text{\tiny def}}{=} 
\frac{1}{2} \int_\mathbb{R} V(x) |q_N(x,a,v,\theta,\mu)|^2\, dx \,.
\]
Although somewhat complicated looking, the equations \eqref{eq:effd}
have a natural interpretation in terms of Hamiltonian
systems: they are the Hamilton-Jacobi equations for the full 
Hamiltonian of \eqref{eq:GP} restricted to the symplectic 
$ 4N$-dimensional manifold of $ N$-solitons -- see \S \ref{sec:eff} 
for details. Of course when $ V \equiv 0 $
the solutions correspond to the exact solutions of \eqref{eq:uxt0}.
The mathematical results of  \cite{holmer-zworski1}, \cite{holmer-zworski}
suggest that the approximation \eqref{eq:uxtV} is valid up to 
a (rescaled) Ehrenfest time \eqref{eq:Ehr}:
\[ \text{\eqref{eq:uxtV} holds for $ 0 < t < C \log ( 1/h )/ h $. } \]
In other words, 
the equations \eqref{eq:effd} give the 
minimal exact effective dynamics valid up to the Ehrenfest time
$ \log(1/h) /h $, where $ h $ is the parameter controlling the 
small variation of the potential, see \eqref{eq:pot}. 
In this work, we show that the approximation errors
$ {\mathcal O} ( h^2 )$ and the Ehrenfest time bound are sharp.
See \cite{KM} for a survey of soliton dynamics under 
integrable systems that have been perturbed. 
\begin{figure}[t]
\centering
\includegraphics[scale=.8]{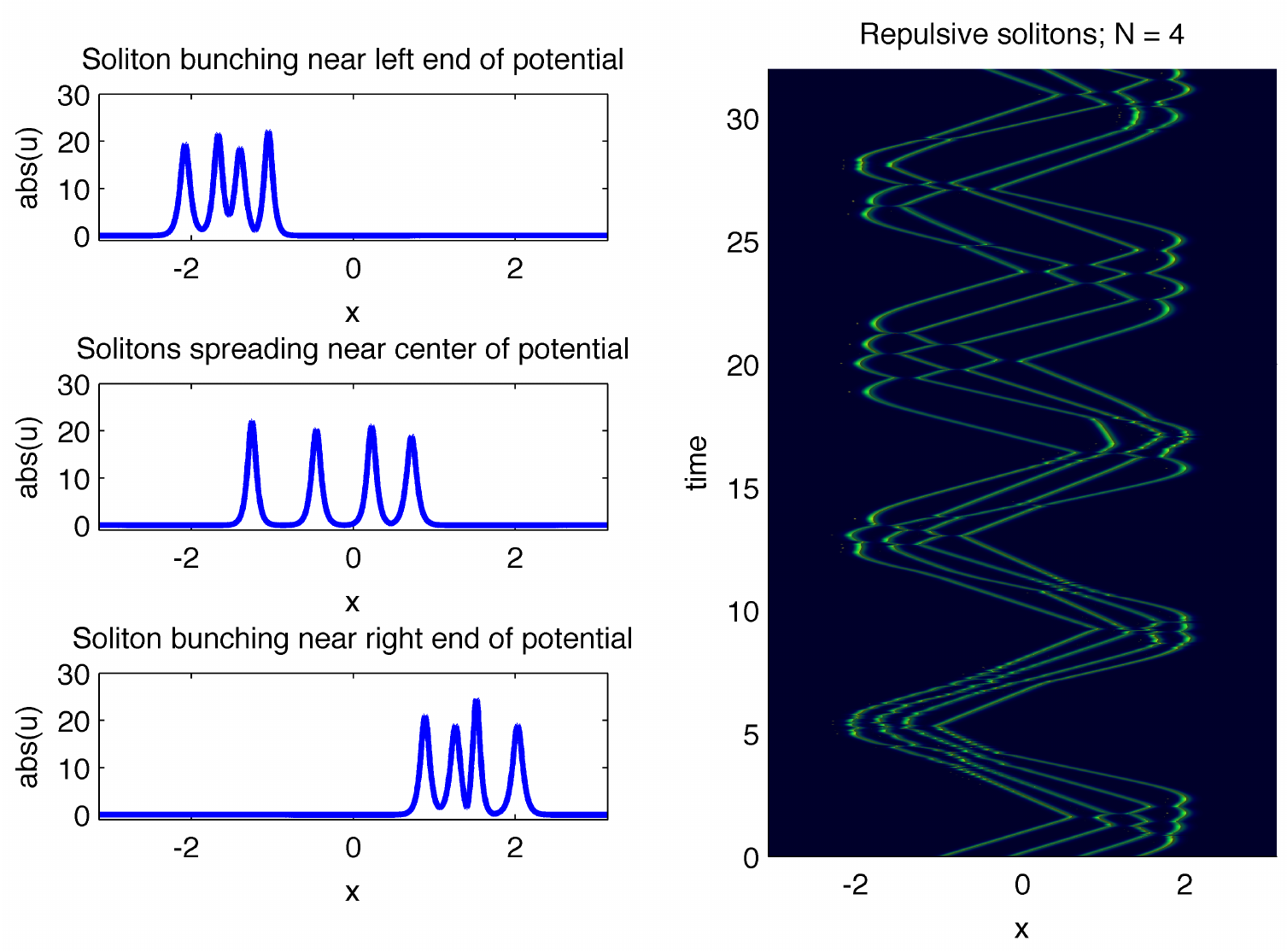}
\caption{Four solitons with alternating phases bunching up and spreading out
in the potential $V(x) =  \left(x/2\right)^6\,.$ 
The full solution is plotted on the right with a bird's eye view. 
The figures on the left are snapshots of that solution.
Due to their alternating phases, the solitons repel 
and never pass through each other.}
\label{fig:bunchspread}
\end{figure}

One motivation for this study is the experimental and theoretical
investigation of soliton trains in Bose-Einstein condensates \cite{strecker}.
We show that the effective dynamics described in \S \ref{sec:eff} is
in qualitative agreement with the behaviour of the matter-wave soliton trains 
-- see Figure \ref{fig:bunchspread}.

The paper is organized as follows: in \S \ref{sec:Ns} we recall the construction of
$ N$-soliton solutions for $ V \equiv 0$ and in \S \ref{sec:eff}, the Hamiltonian
structure of the equation and the derivation of the effective equations of motion.
In \S \ref{effex} we compare the effective dynamics to the behaviour of
solutions to \eqref{eq:GP} and draw some quantitative conclusions. 
Specific potentials similar to those in \cite{strecker} are then discussed in 
\S \ref{sec:appbec}. We investigate effective dynamics for the mKdV equation in
\S \ref{sec:efmk}. Finally, in \S \ref{sec:num} we describe the numerical methods
and compare to other possible approaches.

\begin{figure}[t]
\includegraphics[scale=.7]{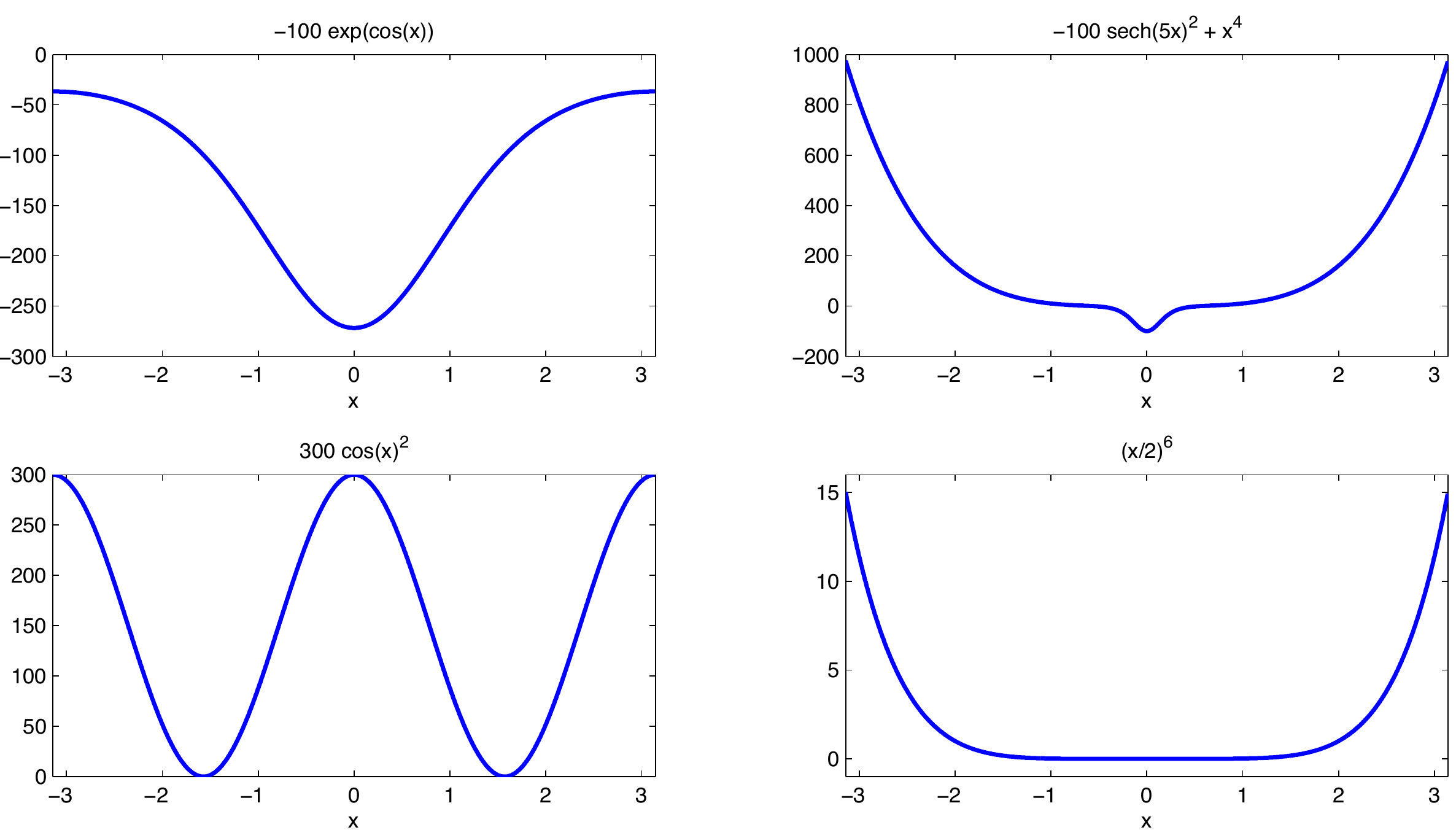}
\caption{A gallery of potentials used for the numerical experiments. 
Since the solitons in the experiments have width approximately $1/10$,
the interesting potentials should have size approximately $100$.
This is suggested by the rescaling (\ref{eq:scalpha}).
The potentials vary on a scale comparable to $1$, 
hence the effective $h$ is approximately $1/10$. 
The exception is the upper right plot where we intentionally
chose a potential which will exhibit some failures of effective dynamics.
In the analysis of errors, for instance in Figure \ref{fig:loglog},
only relative sizes of $h$ matter.}
\label{fig:potplot}
\end{figure}

\section{$N$-solitons for cubic NLS}
\label{sec:Ns}
When $V \equiv 0$, we recover the nonlinear cubic one dimensional 
Schr\"odinger equation, 
which has $N$-soliton solutions with explicit formulas that we
now recall -- see \cite{FT} for a detailed presentation of this
completely integrable equation. 

We will construct functions $q_N(x)$ that depend on $4N$ parameters: 
positions, velocities, phases, and masses:
\begin{equation}
q_N(x) = q_N(x,a,v,\theta,\mu),
\quad a,v,\theta \in \mathbb{R}^{N}, \quad \mu \in (0, \infty)^N.
\end{equation}
Put
\begin{equation} \lambda_j = v_j + i\mu_j\, , 
\ \ \quad \gamma_j(x) = e^{i\lambda_jx}e^{i(\theta_j - v_j a_j)}e^{\mu_j a_j}\,,
\end{equation}
and define matrices 
\[ M(x) \in \mathbb{R}^{N\times N} \,, \ \ 
\ \ \ M_1(x) \in \mathbb{R}^{(N+1)\times (N+1)} \,, \]
by
\begin{equation}
M_{jk}(x) = \frac{1 + \gamma_j(x) \bar{\gamma}_k(x)}{\lambda_j -
 \bar{\lambda}_k} \,, \ \ \ 
M_1 = \left[ \begin{array}{cc} M(x) & \gamma \\ \ & \ \\ 
\vec{1} & 0 \end{array}\right]
\end{equation}
where 
\begin{equation}\gamma = \left[\gamma_1,  \cdots,  \gamma_N \right]^T,
\quad \vec{1} = \left[1, \cdots , 1 \right] \in \mathbb{R}^N.
\end{equation}
Finally,
\begin{equation}
\label{eq:q_N}
q_N(x) \stackrel{\text{\tiny def}}{=}  \frac{\det M_1(x)}{\det M(x) }.
\end{equation}
Remarkably, this gives a solution to \eqref{eq:GP0} with 
$ V \equiv 0 $, the $N$-soliton solution:
\begin{equation}
u(x,t) = q_N(x,a+tv,v,\theta + \frac{t}{2}(\mu^2 + v^2),\mu)\,.
\end{equation}
As one can see from the formula, some restrictions on the parameters apply,
see \cite{FT}.

\section{Effective Dynamics Equations}
\label{sec:eff}

We consider potentials defined on $\mathbb{R}$ that are slowly varying in the sense that
\begin{equation}
\label{eq:pot}  V(x) = W(hx) 
\end{equation}  
where $ W ( x ) $ is $ C^2 $ in $ x $,  and 
\[  | \partial_x^k W ( x ) | \leq C ( 1 + |x| )^N \,,  k \leq 2 \,. \]
where $ C $ and $ N $ are independent of $ h $.
This means that $ h$ is the parameter controlling the slow variation of $ V $.

To obtain an effective dynamics for the evolution we use the
Hamiltonian structure of the equation. In the physics literature
an approach using Lagrangians is more common -- see for instance 
Goodman-Holmes-Weinstein \cite{GHW} and Strecker et al \cite{strecker}.
In the mathematics treatments 
\cite{FrSi},\cite{holmer-zworski},\cite{holmer-zworski1} 
the Hamiltonian approach was found easier to use, 
which we follow here.

The basic claim is that an approximate evolution of $q_N$ 
is obtained by restricting the Hamiltonian flow 
generated by the Gross-Pitaevskii equation to the manifold of $N$-solitons
described in \S \ref{sec:Ns}. 
The Hamiltonian associated with the Gross-Pitaevskii equation is
\begin{equation}
H_V(u)  \stackrel{\text{\tiny def}}{=} 
\frac{1}{4}\int(|\partial_xu|^2 - |u|^4)\, dx
+ \frac{1}{2} \int V|u|^2
\end{equation}
with respect to the symplectic form 
\begin{equation}
\omega(u,v) = \mbox{Im}  \int u\bar{v}\,.
\end{equation}
The manifold of solitons, $M_N$, is $4N$-dimensional 
and equipped with the restricted symplectic 
form given by the sum of forms for single solitons:
\begin{equation}
\label{eq:omegaM}
\omega_M \stackrel{\text{\tiny def}}{=} \omega|_M = 
\sum_{j=1}^N (\mu_j dv_j \wedge da_j + v_j d\mu_j \wedge da_j + 
d\theta_j \wedge d\mu_j)\,.
\end{equation}
$H_V$ restricted to $M_N$ is 
\begin{equation}
H_N \stackrel{\text{\tiny def}}{=} H_V|_{M_N}(a,v,\theta,\mu) = \sum_{j=1}^N \left( \frac{\mu_j v_j^2}{2} - \frac{\mu_j^3}{6} \right)
+ V_N(a,v,\theta,\mu)\,, 
\end{equation}

\begin{equation}\label{eq:V-N}
\mbox{where} \quad V_N(a,v,\theta,\mu) \stackrel{\text{\tiny def}}{=} 
\frac{1}{2} \int_\mathbb{R} V(x) |q_N(x,a,v,\theta,\mu)|^2\, dx\,.
\end{equation}
The effective dynamics is given by the flow of the Hamilton vector field
of $H_N$ on the manifold $M_N$. 
That vector field, $\Xi_{H_N}$, is defined using the symplectic form
(\ref{eq:omegaM}):
\begin{equation}
dH_N = \omega_M(\cdot, \Xi_{H_N})\,.
\end{equation} 

A computation based on this gives the following ordinary differential equation 
for the parameters $a, v, \theta$ and $\mu$,  
called the effective dynamics:
\begin{equation}\label{eq:effdyn} 
\begin{split}
\dot{v}_j & = -\mu_j^{-1} \partial_{a_j} H_N - \mu_j^{-1}v_j 
\partial_{\theta_j} H_N = -\mu_j^{-1} (\partial_{a_j} V_N + v_j \partial_{\theta_j} V_N)\,, 
\\ 
\dot{a}_j & = \mu_j^{-1} \partial_{v_j} H_N 
=  v_j + \mu_j^{-1} \partial_{v_j}V_N \,, \\
\dot{\mu}_j & = \partial_{\theta_j} H_N = \partial_{\theta_j} V_N\,,  \\
\dot{\theta}_j & -= \mu_j^{-1} v_j \partial_{v_j} H_N - \partial_{\mu_j}H_N 
=   v_j^2/2 + \mu^2_j/2 + \mu_j^{-1} v_j \partial_{v_J} V_N - \partial_{\mu_j} V_N \,. 
\end{split}
\end{equation}

We remark that one can scale the Gross-Pitaevskii equation (\ref{eq:GP}) 
in the following way: Consider any function $u(x,t)$, scaling parameter $\alpha$,
let $\tilde x = \alpha x\,, \tilde t = \alpha^2 t$, and define the new function 
\begin{equation}\label{eq:scalpha}
\tilde u(\tilde x,\tilde t)  \stackrel{\text{\tiny def}}{=} \frac{1}{\alpha}u(x,t)\,.
\end{equation}
Then if $u(x,t)$ satisfies (\ref{eq:GP}) with the potential $V(x)$,
$\tilde u(\tilde x, \tilde t)$ also satisfies (\ref{eq:GP}) with 
the new potential 
\begin{equation}
\label{scaleV}
\tilde V(\tilde x) = \frac{1}{\alpha^2} V\left(\frac{ \tilde x}{\alpha}\right) \,. 
\end{equation}
This means that if we deal with a soliton of width comparable with $\alpha$,
the potentials for which interesting dynamics should appear should have size 
approximately $\alpha^{-2}$ and the slowly varying factor replaced by $h/\alpha$.

The effective dynamics equations (\ref{eq:effdyn}) scale similarly:
if $\left( a(t), v(t), \theta(t), \mu(t) \right)$ satisfies (\ref{eq:effdyn}) 
and we define $\tilde x, \tilde t,$ and $\tilde V$ as above, then 
$\left( \tilde a(\tilde t), \tilde v(\tilde t), \tilde \theta(\tilde t), \tilde \mu(\tilde t) \right)$
also satisfies (\ref{eq:effdyn}) with  
\begin{equation}\label{eq:effdynsca} 
\tilde a(\tilde t) = \alpha a(t)\,, \quad \tilde v(\tilde t) = \frac{v(t)}{\alpha}\,, \quad
\tilde \theta(\tilde t) = \theta(t)\,, \quad
\tilde \mu(\tilde t) = \frac{\mu(t)}{\alpha}\,.
\end{equation}

The scalings (\ref{eq:scalpha}) and (\ref{eq:effdynsca}) are related in the following way:
if
 $$u(x,t) = q_N(x,a(t),v(t),\theta(t),\mu(t))\,,$$ 
then 
$$\tilde u(\tilde x, \tilde t) = q_N(\tilde x,a(\tilde t), v(\tilde t), \theta(\tilde t),
\mu(\tilde t))/\alpha = 
q_N(x, \tilde a(\tilde t), \tilde v(\tilde t),\tilde\theta(\tilde t), \tilde \mu(\tilde t))\,.$$
\begin{figure}[t]
\begin{tabular}{c c}
\includegraphics[scale=.55]{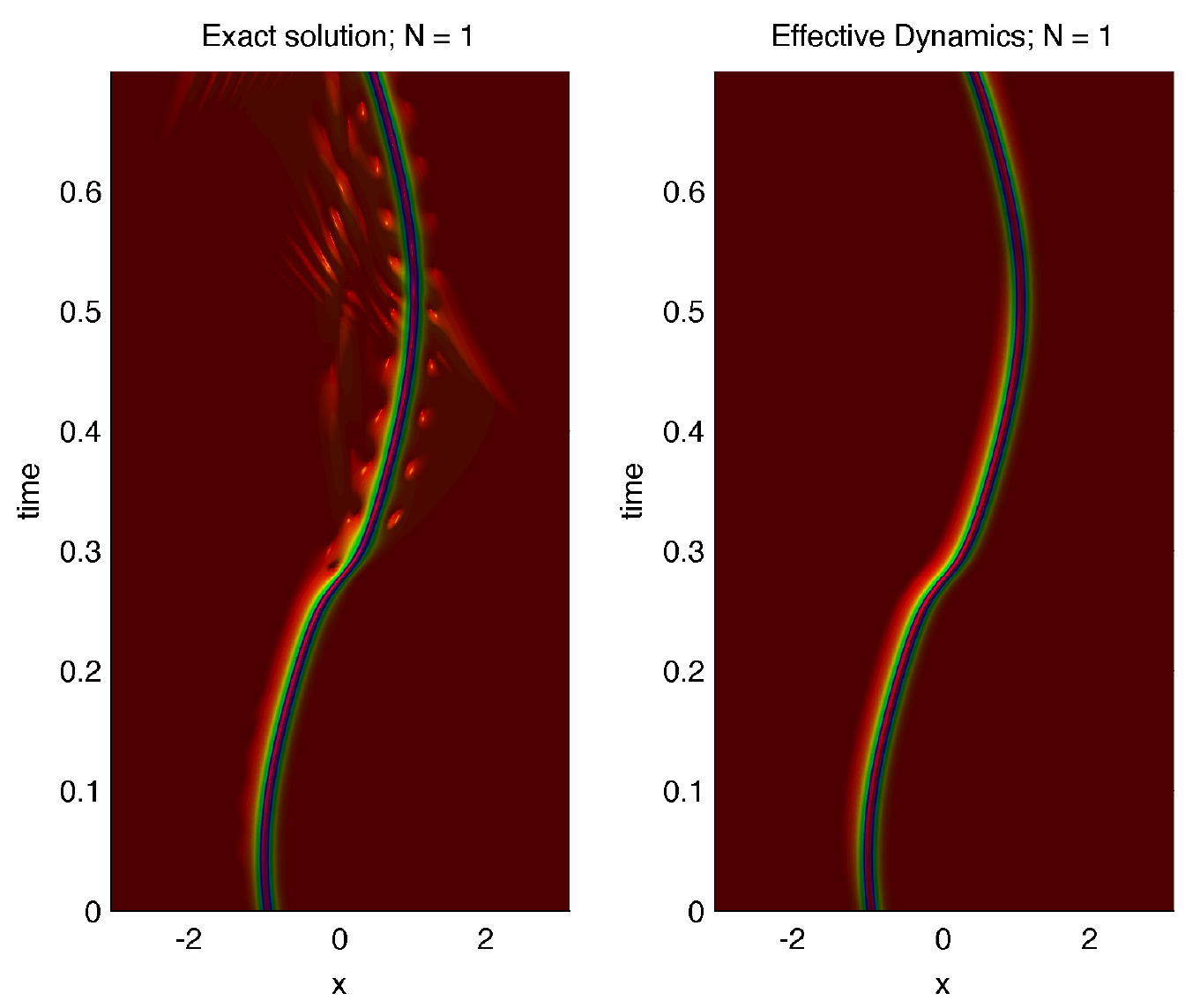} & \includegraphics[scale=.55]{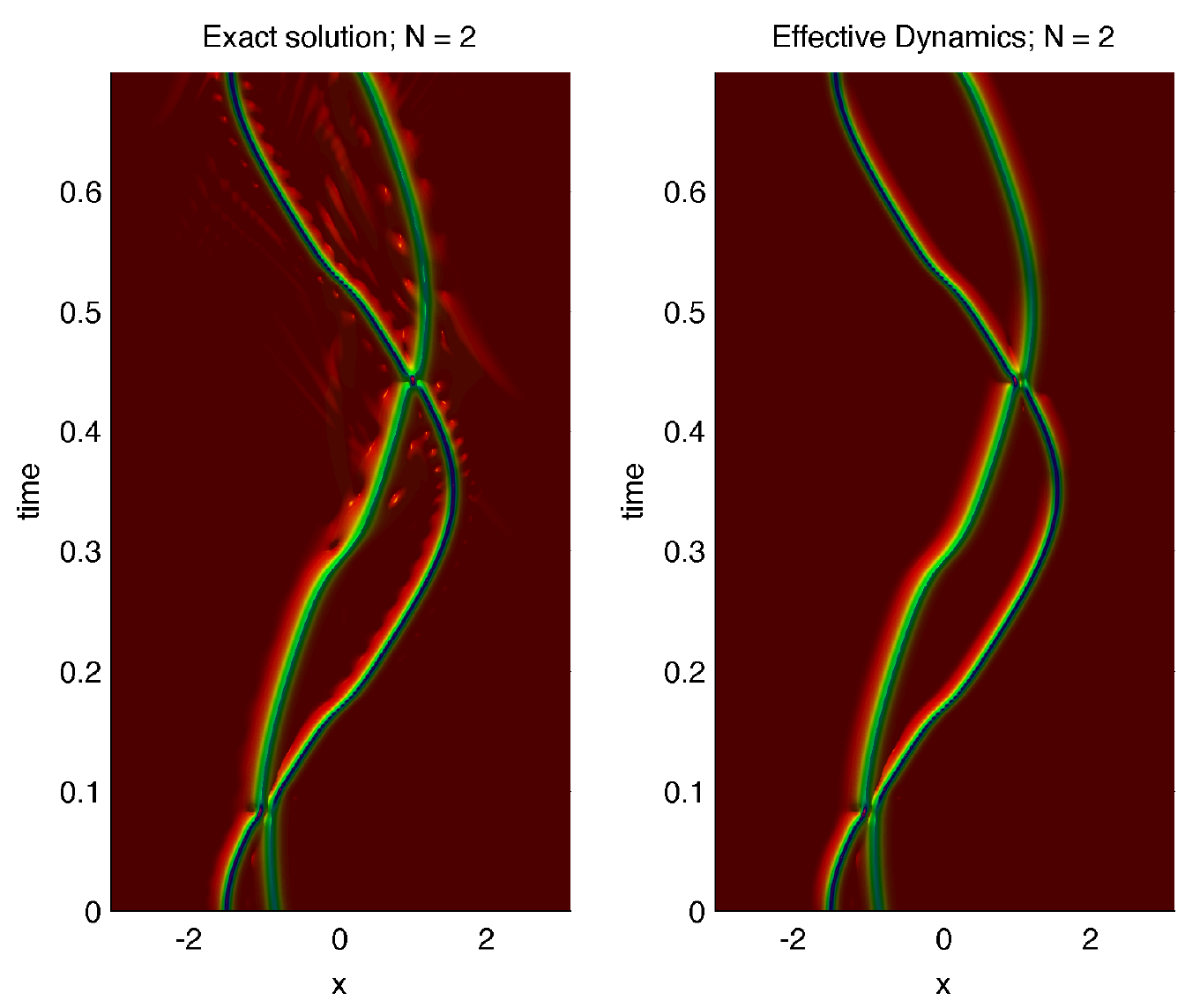} \\
\includegraphics[scale=.55]{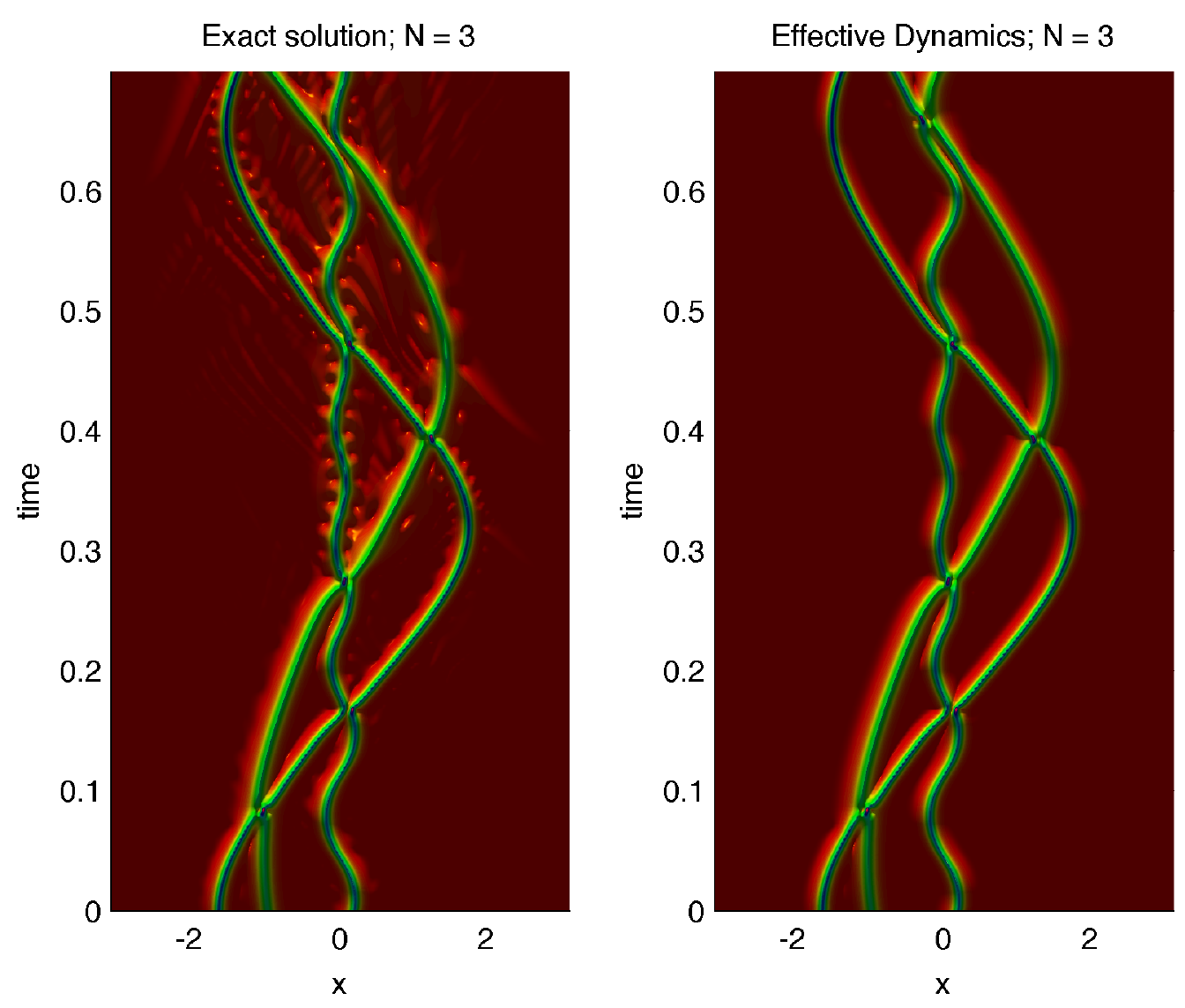} & \includegraphics[scale=.55]{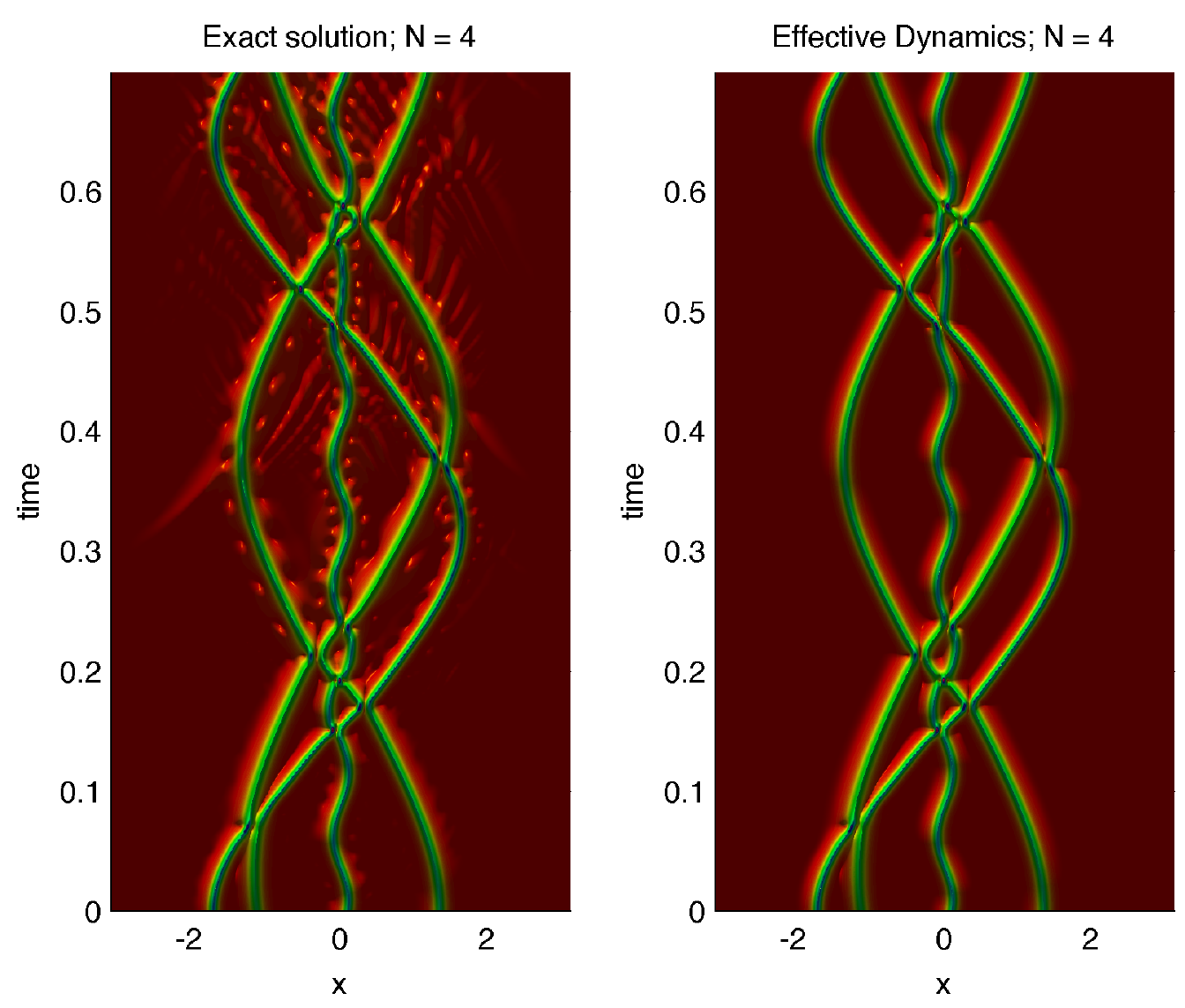} 
\end{tabular}
\caption{A side-by-side comparison of the absolute value of the
exact solution of (\ref{eq:GP}) versus the
effective dynamics (\ref{eq:effdyn}) for 1,2,3, and 4 solitons with potential
$W(x) = -100\, \text{sech}^2 (5x) + 10x^4 \, .$ 
The sharpness of the $\text{sech}^2(5x)$ term creates
clearly visible discrepancy between the two solutions.}
\label{fig:comp}
\end{figure}

\section{Comparison of effective and exact dynamics}
\label{effex}
For given values $a_0, v_0, \theta_0, \mu_0$ in $\mathbb{R}^N$, 
we consider the solution $u(\cdot,t)$ of (\ref{eq:GP}) with initial data 
$q_N(\cdot,a_0,v_0,\theta_0,\mu_0)$ and the solutions $a(t),v(t),\theta(t),\mu(t)$ of 
the effective dynamics equations (\ref{eq:effdyn}) with initial values $a_0,v_0,\theta_0,\mu_0$.
In the following discussions we will refer to $u(\cdot,t)$ as the exact solution and
$q_N(\cdot,a(t),v(t),\theta(t),\mu(t))$ as the effective dynamics.

Holmer and Zworski \cite{holmer-zworski} proved that in the case $N=1$,
\begin{equation}
\label{eq:apprGP} \|  u ( \cdot,t ) - q_N ( \cdot, a(t),v(t),\theta(t),\mu(t) ) \|_{H^1}  = 
Ch^{2-\delta}\,, \mbox{ for } t < \frac{\delta \log(1/h)}{Ch}\,,
\end{equation}     
where $\delta \in (0,1/2)$ can be chosen, and 
where $C$ depends only on the potential and initial velocity of the soliton,
but not on $\delta$.
The $H^1$ norm measures the size of the function and its spatial derivative
in $L^2$:
\[  \| v \|_{H^1}^2 \defeq \|v\|_{L^2}^2 + \|\partial_xv\|_{L^2}^2\,, 
\qquad  \| v \|_{L^2}^2  \defeq  \int_\RR |v ( x )|^2 dx \,. \]    
This norm measures the energy of the solution. 

The limiting time $ \log(1/h)/h $ is the Ehrenfest time discussed
int \S \ref{sec:intro}.

It is expected that this result also holds for $N>1$.
This is suggested by \cite{holmer-zworski1}, which proves the analagous theorem
for the modified Korteweg-de Vries (mKdV) equations for case $N=1,2$, see \S \ref{sec:efmk} below.
However, the methods of \cite{holmer-zworski1} do not fully apply to the 
case of the Gross-Pitaevsky equation \eqref{eq:GP0}. 
Also, even in the case of mKdV and $ N = 2$, multiple soliton interactions are
not theoretically understood. 
All these considerations provided a strong motivation 
for this numerical study.

We note that for $N>1$, it was conjectured in \cite{holmer-zworski1} 
that the error bounds (\ref{eq:apprGP}) will hold 
not only in the $H^1$ norm, but for the $H^N$ norm, 
which measures the the size of a function and its first $N$ derivatives,
where $N$ is the number of solitons.
This has been proven in the mKdV case with $N=2$ \cite{holmer-zworski1}.
For our numerical experiments, we consider only the $H^1$ norm. 

We present numerical simulations to show that the result (\ref{eq:apprGP}) holds for $N>1$ 
in the following three sections:  
In \S\ref{subsec:ancs} we choose initial data and two potentials that demonstrate the power
and limitations of the effective dynamics equations, regardless of the number of solitons. 
Using the same initial data and one of the potentials from \S\ref{subsec:ancs},
in \S\ref{subsec:htozero} we verify that the $\mathcal{O}(h^{2-\delta})$ error estimate in 
(\ref{eq:apprGP}) holds as $h \to 0$ for a fixed time interval. 
We then turn to the $\log(1/h)/h$ timescale, or Ehrenfest timescale,
in \S\ref{subsec:et} to show that it is the 
appropriate timescale for which we can expect (\ref{eq:apprGP}) to hold.

\subsection{A numerical case study}
\label{subsec:ancs}
We consider initial data 
$q_N(\cdot, \bar{a}_N, \bar{v}_N, \bar{\theta}_N, \bar{\mu}_N)$,
where $\bar{a}_N = (a_1,\dots, a_N)\, \ , N = 1,2,3,4$ and 
$\bar{v}_N, \bar{\theta}_N,$ and $\bar{\mu}_N$ are similarly defined with
\begin{equation}
\begin{split}
(a_1,a_2,a_3,a_4) & = (-1, -1.5, 0, 1) \\
(v_1,v_2,v_3,v_4) & = (-2, 0, 3, 0) \\ 
(\theta_1, \theta_2, \theta_3, \theta_4) & = (\pi/3, 0, -3, -5) \\
(\mu_1, \mu_2, \mu_3, \mu_4) &  = (17, 25, 23, 19) \\
\end{split}
\end{equation}

The positions and masses are chosen to satisfy a numerical requirement that
$q_N(\cdot, \bar{a}_N, \bar{v}_N, \bar{\theta}_N, \bar{\mu}_N)$ is close to 0
outside of $(-\pi, \pi)$, our numerical domain (see \S \ref{sec:num}). 
Rescaling the solution as in (\ref{eq:scalpha}) or enlarging the numerical domain
allows for data that does not satisfy this numerical requirement.
The initial data is otherwise chosen arbitrarily. 

We first consider the potential 
$$V_1(x) = -100\,e^{\cos x}\,,$$
see Figure \ref{fig:potplot}. 
The factor $-100$ is chosen to create a deep enough well so that the 
solutions remain in $(-\pi, \pi)$, 
but the potential is otherwise chosen arbitrarily. 

We compute the exact solution and the effective dynamics solution for 
$N = 1,2,3,4$ up to time $t = 1$, which is chosen to allow for multiple soliton
interactions. We plot the solution for 4 solitons in Figure \ref{fig:perfectmatch},
where we observe very little difference between the exact 
and effective dynamics solutions. 
An equally small amount of discrepancy between the solutions
was observed for $N = 1,2,3$.

Next, we consider the same initial data as above with the
potential $$V_2(x) = -100\, \text{sech}^2 (5x) + 10\,x^4 \, ,$$ 
see Figure \ref{fig:potplot}. 
This potential is chosen to be outside of the slowly varying regime for which 
the effective dynamics give good approximations. 
This is due to the $-100\text{sech}^2(5x)$ term, 
which creates a sharp dip roughly the width of the 
solitons we are studying; 
thus we do not expect the exact solution to maintain its soliton structure very well. 
This causes the clearly visible discrepancies 
between the effective dynamics and the exact solution in Figure \ref{fig:comp}.
The $10x^4$ term ensures the solutions remain on the interval $(-\pi, \pi)$.

We compute the solutions for $N = 1,2,3,4$ up to time $t = 0.7$, 
which is again chosen to allow for multiple soliton interactions.
Figure \ref{fig:comp} displays these solutions and
demonstrates that the effective dynamics captures the true motion, 
regardless of the number of solitons and regardless of multiple soliton interactions. 
In the experiments presented in this paper, we only consider $N \le 4$, 
but we have observed good agreement between the exact solution and the effective dynamics
for $N \le 7$. 
We did not investigate futher due to increasing computational time needed to solve
(\ref{eq:effdyn}).

\begin{figure}[t]
\begin{tabular}{c c}
\includegraphics[scale=.65]{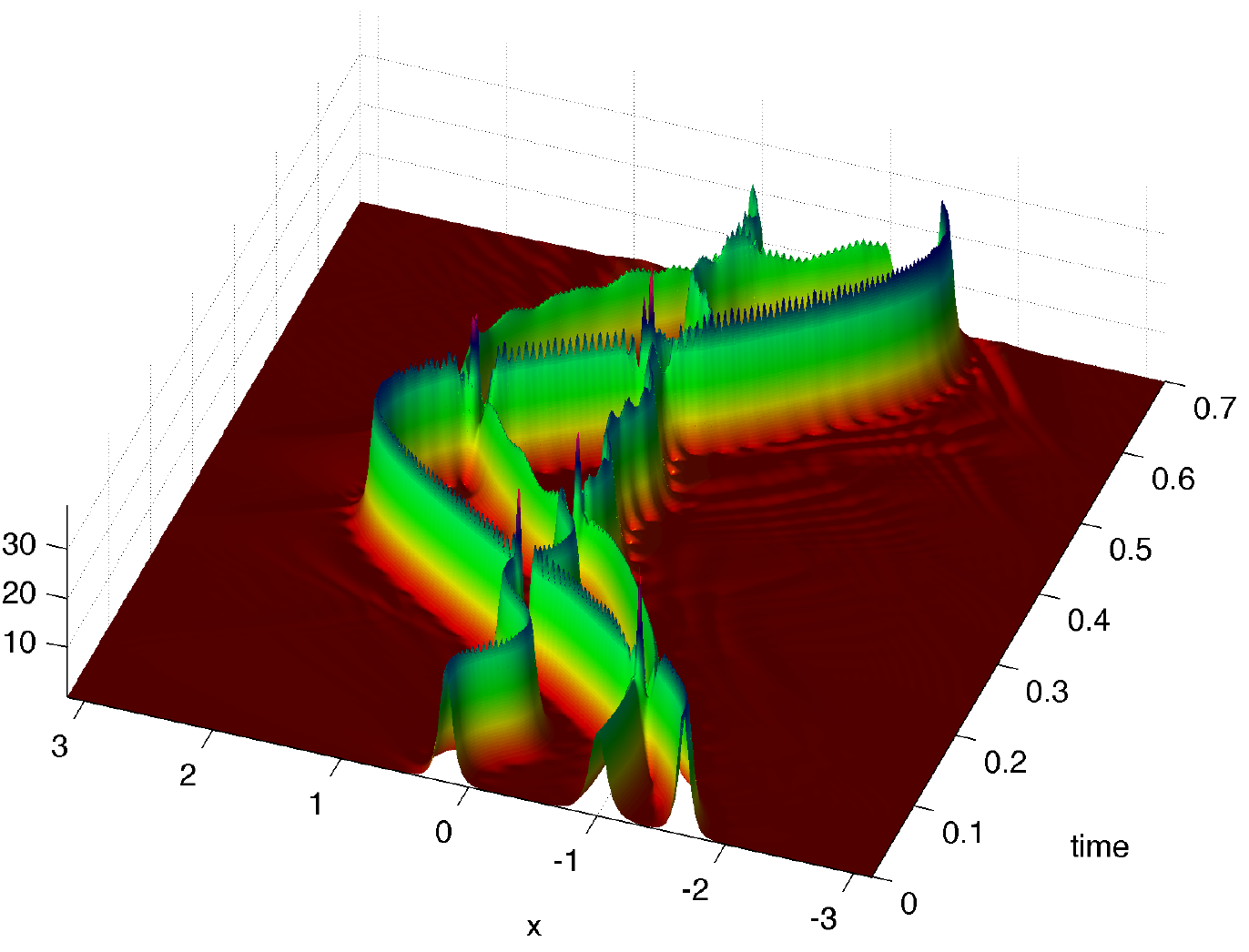} & \includegraphics[scale=.53]{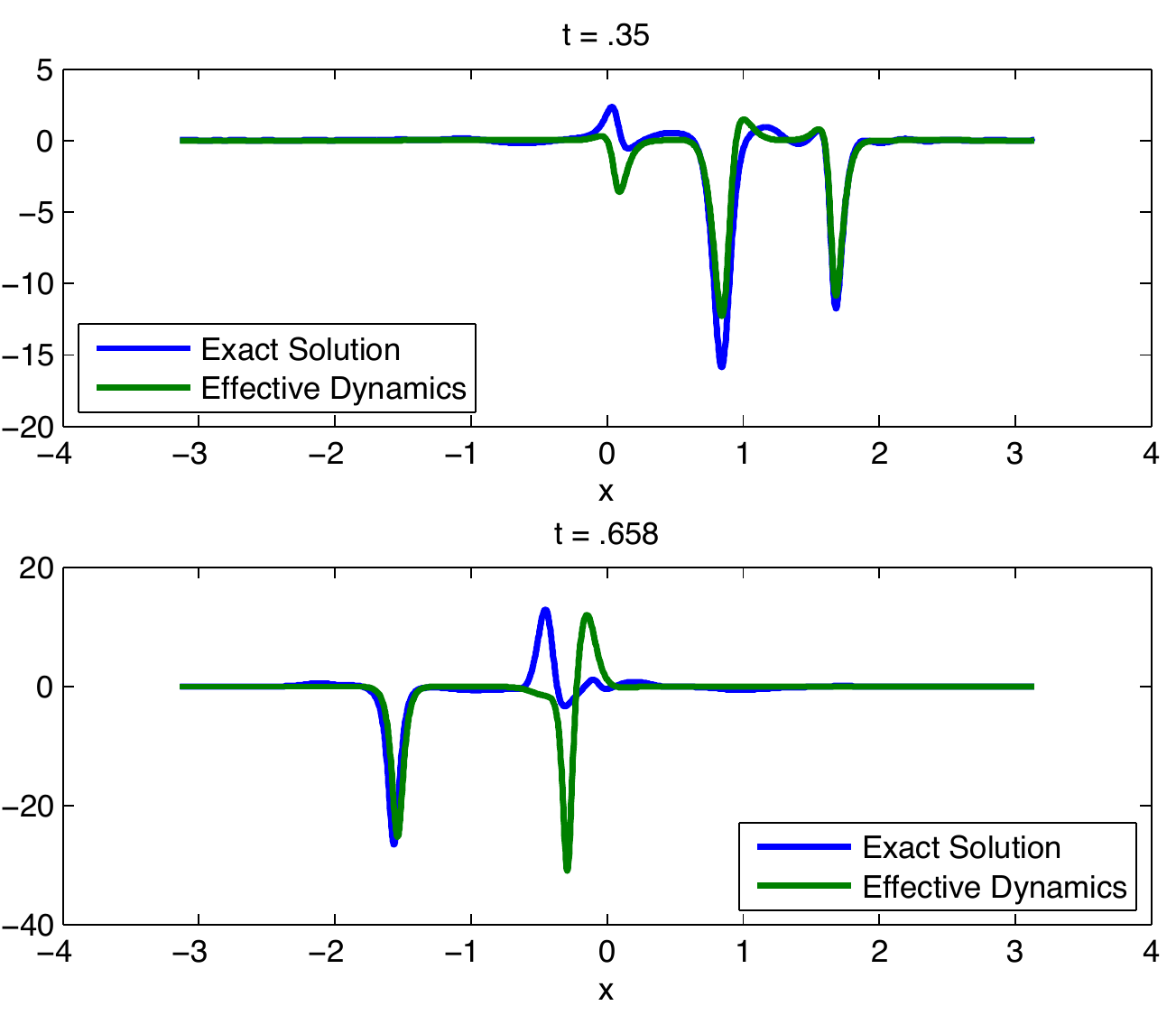}
\end{tabular}
\caption{The plot on the left is a different view of the exact solution with 3 solitons 
shown in Figure \ref{fig:comp}. The plot on the right compares the real parts of the 
$N = 3$ exact solution with the effective dynamics solution at times $t = 0.35$ and $t = 0.658$.}
\label{fig:keepshape}
\end{figure}

We note that for $N \ge 2$ the phases of the solitons 
are crucial in determining the interaction between solitons.
In Figure \ref{fig:comp}, we see that for $N = 3$,
at approximately $t \approx 0.65$,
two solitons that appear to bounce off each other in the exact solution
instead appear to cross in the effective dynamics.
This discrepancy seems to be due to differences between exact phases
and effective phases.
In Figure \ref{fig:keepshape} we are able to see large deviation in the phases between
the exact solution and effective dynamics by comparing the real part of the solutions.

\subsection{Quantitative study of the error as $h \to 0$}
\label{subsec:htozero}
We investigate the $\mathcal{O}(h^{2-\delta})$ error between the exact solution
and the effective dynamics on a fixed time interval. 
The estimate that gives rise to (\ref{eq:apprGP}) is 
\begin{equation}
\label{eq:fullerr} 
\|  u ( \cdot,t ) - q_N ( \cdot, a(t),v(t),\theta(t),\mu(t) ) \|_{H^1}  \le
Ch^2e^{Cht}\,.
\end{equation}     
If $t = \delta \log(1/h)/(Ch)$, then the RHS reduces to $Ch^{2-\delta}$.
When dealing with fixed length of time or even time of size 
$\mathcal{O}(1/h)$ the RHS is $\mathcal{O}(h^2)$ 
and that form of error will be shown to be optimal.

\begin{figure}[t]
\centering
\includegraphics[scale=.5]{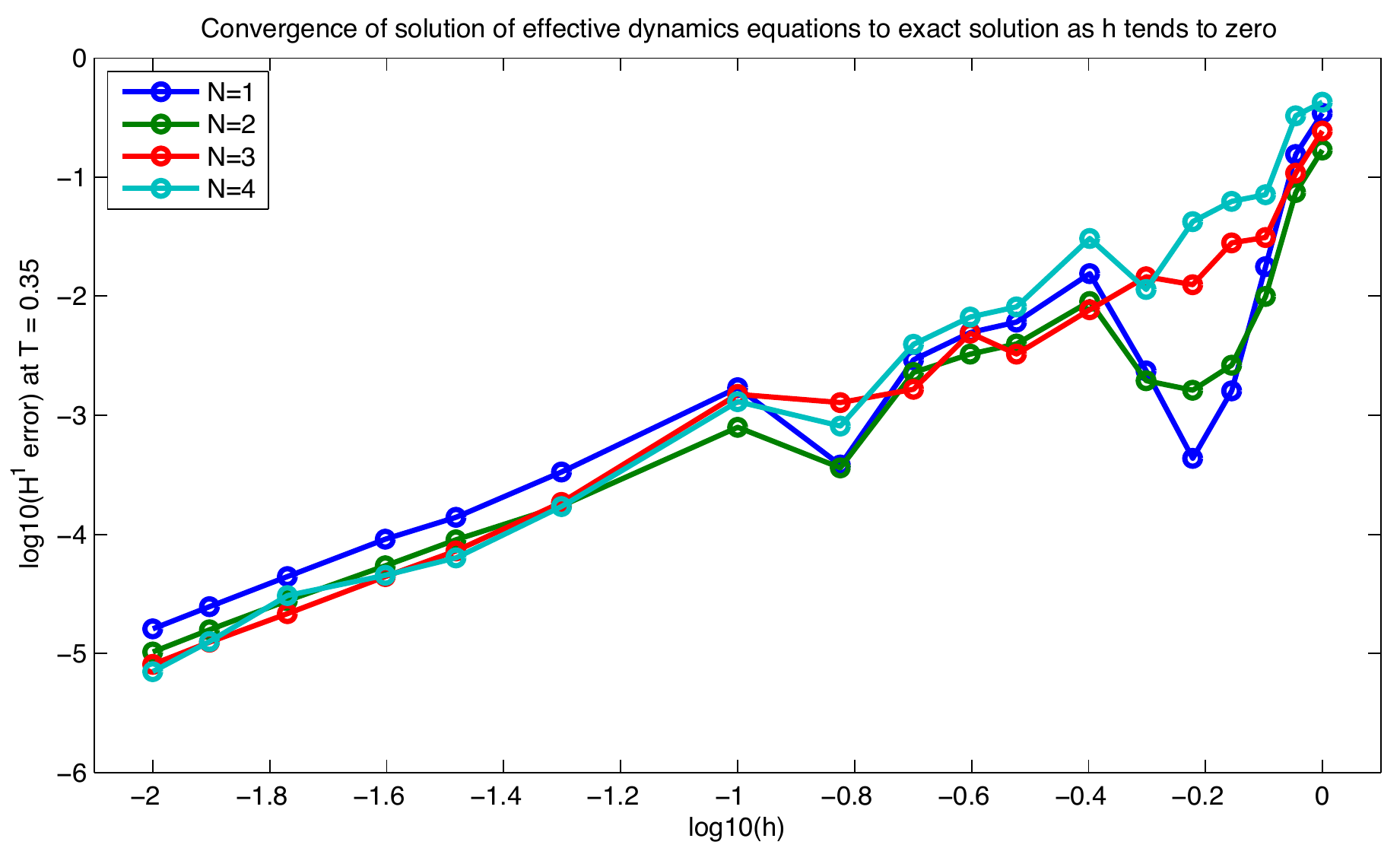}
\caption{A log-log plot of the $H^1$ error, relative to
the $H^1$ norm of the initial data, between the exact solution of \ref{eq:GP} and 
the $N$-soliton evolving according to the effective dynamics equations \ref{eq:effdyn},
 as function of $h$. Here, the potential is $V(x) = W(hx)$,
where $W(x) = -100\, \text{sech}^2 (5x) + 10x^4\,$.
For smaller values of $h$, the slope of the lines approaches 2, 
in agreement with the theoretical upper bound on the error in (\ref{eq:uxtV}). }
\label{fig:loglog}
\end{figure}

We reconsider our second potential from \S \ref{subsec:ancs}, 
but add the slowly varying parameter, $h$:
$$V(x) = W(hx)\,, \quad W(x) = -100\, \text{sech}^2 (5x) + 10x^4\, ,$$ 
and explore the $H^1$ error, 
relative to the $H^1$ norm of the initial data,
 between the exact solution and effective dynamics
as $h \to 0$.
We expect that as $h$ becomes small enough, 
the equation will enter the slowly varying regime and
display $\mathcal{O}(h^2)$ error.
Indeed, the log-log plot in Figure \ref{fig:loglog} demonstrates the error is bounded by
$C_N h^2$ as $h \to 0$, where the constant $C_N$ varies only slightly
between different values of $N$. 

We fit the data from Figure \ref{fig:loglog} 
to a line using the 6 smallest values of $h$.
\begin{center}
\begin{tabular}{|c|c|c|c|c| }
\hline
N  & 1 & 2 & 3 & 4 \\
\hline
Slope & 1.86 & 1.76 & 1.91 & 1.86 \\
\hline
$C_N$ & -1.07 & -1.46 & -1.27 & -1.37 \\
\hline
\end{tabular}
\end{center}
Thus, we conclude that 
the error is approximately $\mathcal{O}(h^2)$.

\subsection{Ehrenfest time}
\label{subsec:et}
We now investigate the length of time for which 
the effective dynamics approximation is accurate.
In (\ref{eq:fullerr}), we recalled that the error is bounded by
$Ch^2 e^{Cht}$
and hence the approximation breaks down
at the Ehrenfest time, $$T(h) \sim \log ( 1/h )/ h\, .$$
We have already verified for small $h$ and fixed time
this error behaves as $\mathcal{O}(h^2)$. 
Thus we focus on observing exponential growth in the error as a function of time,
and verifying that it is of the form $\mathcal{O}(e^{Cht})$. 

For this we must choose a potential and initial data to exhibit exponential
instability.
We are motivated by Newton's equations for $V(a) = -a^2/2 $: 
\begin{equation}
\label{eq:newton}
\dot a = v\, ,\quad \dot v = a\, , 
\quad v = v_0\cosh t + a_0 \sinh t\, , \quad a = a_0 \cosh t + v_0 \sinh t \, .
\end{equation}
In this case, we have exponential instability of classical dynamics.
This suggests choosing potentials with a non-degenerate maximum and working near
the unstable equilibrium points. 

\begin{figure}[t]
\includegraphics[scale=.6]{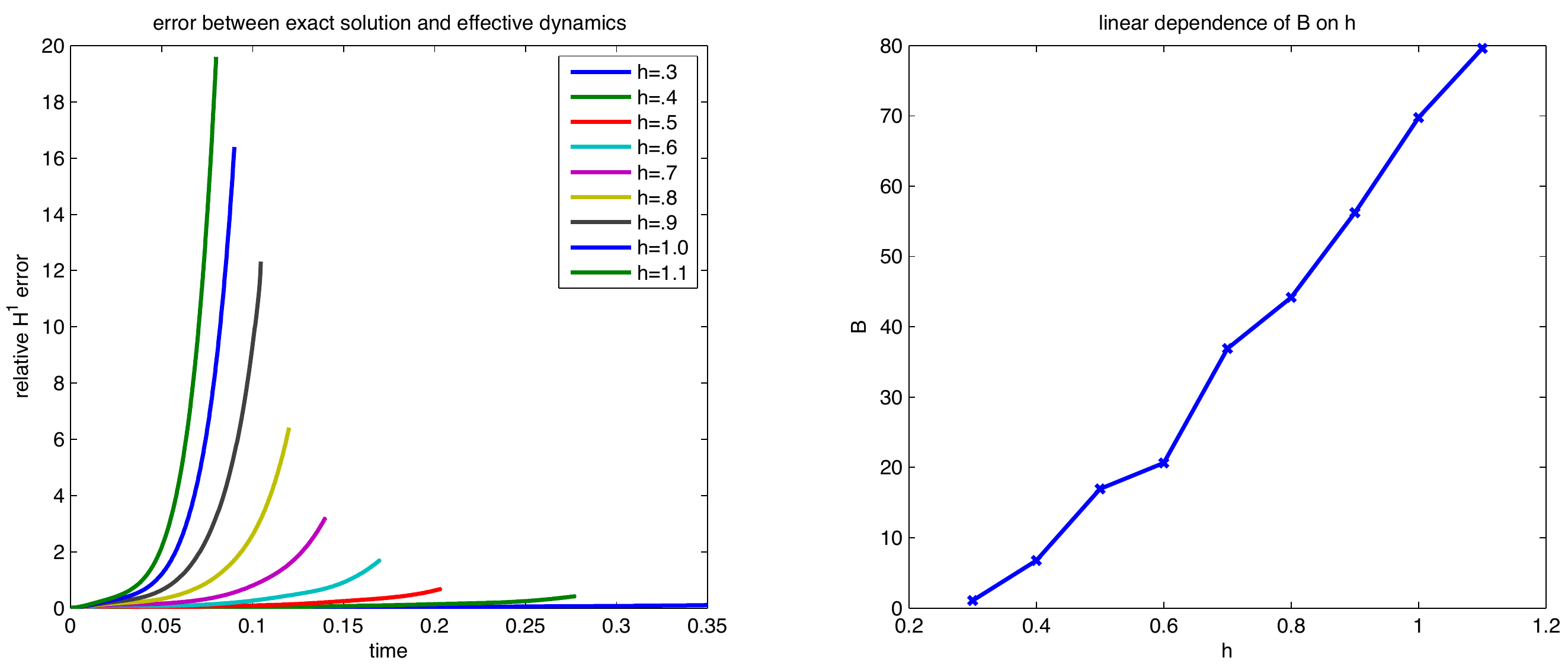}
\caption{The plot on the left shows $H^1$ error,
relative to the $H^1$ norm of the initial data,
between the exact solution of (\ref{eq:GP}) and the effective dynamics 
for a single soliton sliding down the concave potential 
$V(x) = W(hx)\,, \mbox{ where } W(x) = 1000\, x^2\, .$ 
The error is plotted as a function of time and for several values of the parameter $h$. 
On the right sight,
$B$ is plotted as a function of the parameter $h$,
when the errors from the plot on the left are fitted to a curve of the form
$A(e^{Bt} + C)$. We expect $B$ to depend linearly on $h$.}
\label{fig:expdiv}
\end{figure}

Hence we will investigate solutions to (\ref{eq:GP}) with potential
$$V(x)= W(hx)\,, \mbox{ where } W(x) = -1000\,x^2\, .$$ 
Figure \ref{fig:expdiv} below demonstrates
exponential divergence between the exact solution to
(\ref{eq:GP}) and the effective dynamics
for several values of $h$ and a single soliton initial condition $q_1(x,.1,0,0,15)$.

\begin{figure}[t]
\includegraphics[scale=.7]{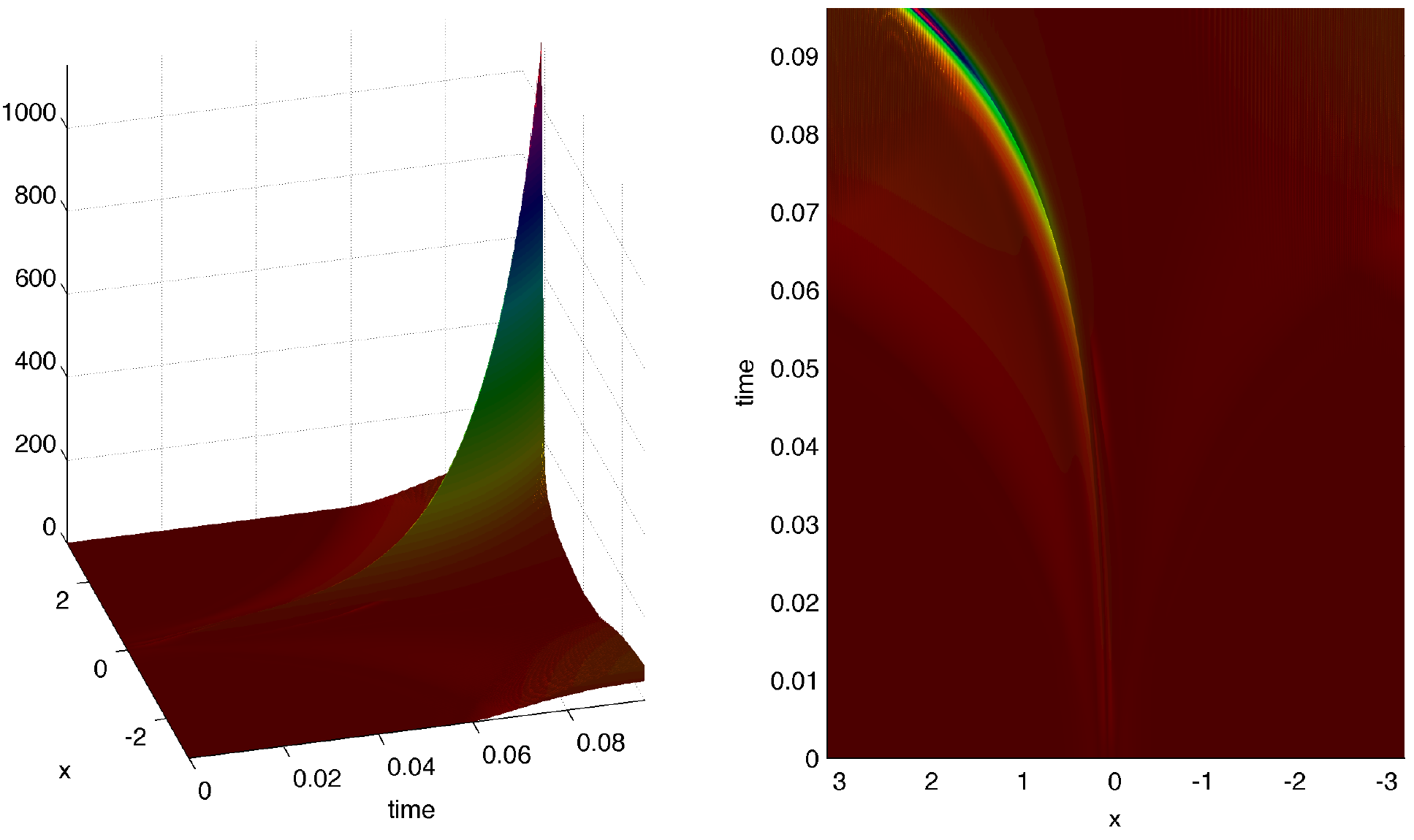}
\caption{Surface plots of the error 
between the exact solution of (\ref{eq:GP}) and the effective dynamics 
for a single soliton sliding down the concave portion of the potential 
$V(x) = W(hx)\,, \mbox{ where } W(x) = 1000\, x^2\, ,$ 
as in Figure \ref{fig:expdiv}. 
We have plotted the absolute value of the difference between the
spatial derivatives between the two solutions.\label{fig:expdiv3d}}
\end{figure}
We fit the plots shown in Figure \ref{fig:expdiv}
to a function of the form $A(e^{Bt} + C)$, 
for the time period when the soliton's position was between $x=.15/h$ 
and $x=1.2/h$.
This range was observed to be a region where exponential increase dominated 
the error and before the soliton approached the numerical boundary.

In Figure \ref{fig:expdiv} we observe a linear dependence of $B$ on $h$,
in agreement with (\ref{eq:fullerr}). 
This indicates that for certain potentials
the Ehrenfest time $C \log ( 1/h )/ h$ is 
the appropriate bound for the length time we expect the effective 
dynamics to give a good approximation to (\ref{eq:GP}).
We note that in our experiments with other potentials we often
observe a linear increase in error which would 
correspond to a timescale of $C/h^2$ instead of the 
Ehrenfest time $C \log( 1/h)/h.$

\vspace{2 mm}

\section{Application to Bose-Einstein Condensates}
\label{sec:appbec}
Strecker, Partridge, Truscott and Hulet  \cite{strecker}
discovered that Bose-Einstein condensates form stable soliton trains
while confined to one-dimensional motion. 
When set into motion in a suitably chosen optical trap, 
a Bose-Einstein condensate forms multiple soliton formations which 
exist for multiple oscillatory cycles without being destroyed
by dispersion or diffraction.
We can observe this same behavior numerically, 
using the effective dynamics equations. 
We choose a potential of the type described in \cite{strecker}
$$V(x) = \left(\frac{x}{2}\right)^6$$
and set $N=4$, which was the most frequent case in their experiment. 
Strecker \emph{et al} inferred that the repulsive behavior of the solitons
indicated alternating phases.
Their argument was based on considering a certain reduced Lagrangian.

In our numerical experiment we put $\bar{\theta} = (0, \pi ,0 ,\pi)$
and then set the four solitons in motion with the same velocity near
the center of the potential. 
Similar to \cite{strecker}, we observe bunching and spreading
of the soliton train for several oscillations. See Figure \ref{fig:bunchspread}.

\vspace{2 mm}

\section{Effective dynamics  for the mKdV equation}
\label{sec:efmk}

The mKdV equation 
\begin{equation}
\label{eq:mkdv0}
\partial_tu = -\partial_x(\partial_x^2 u + 2u^3) \, ,
\end{equation}
like the nonlinear Schr\"odinger equation, 
has soliton solutions and a Hamiltonian structure. 
Holmer, Perelman, and Zworski \cite{holmer-zworski1} 
derived effective dynamics equations
for the mKdV equation with a slowly varying potential
\begin{equation}
\label{eq:mkdv}
\partial_tu = -\partial_x(\partial_x^2 u - b(x,t)u + 2u^3),\quad  b(x,t) = b_0(hx,ht)
\end{equation}
and proved a result analogous to the (\ref{eq:apprGP}) for $N = 1,2$:
the $H^N$ error between the solution of (\ref{eq:mkdv}) and 
its associated effective dynamics with $N$-soliton initial data is bounded by
\begin{equation}
\label{eq:errbdsmkdv}
Ch^2 e^{Cht} \,, \ \ \ \ \ 
\text{ for } \ \  t < \frac{C}{h}  \log \frac{1}{h},
\end{equation}

\begin{figure}[t]
\includegraphics[scale=.65]{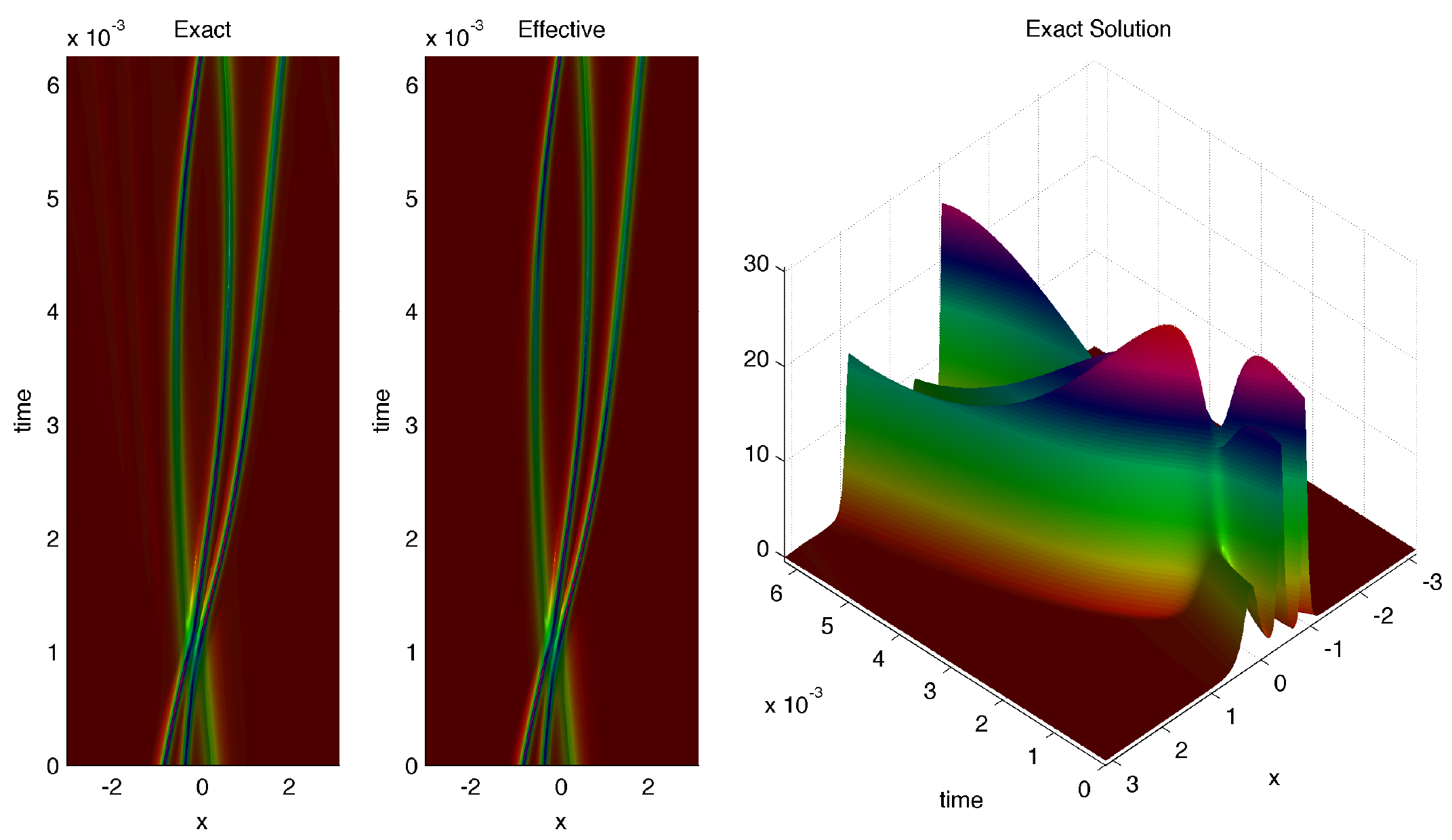}
\caption{The left plot shows a side-by-side 
comparison of the exact solution of the mKdV equation (\ref{eq:mkdv})
and the effective dynamics solution for 3 solitons with potential
$b(x) = 300\cos^2 x$.
No discrepancy between the two solutions is visible.
The right plot displays the exact solution from a different angle.}
\label{fig:birdseyemkdv}
\end{figure}

\begin{figure}[t]
\includegraphics[scale=.7]{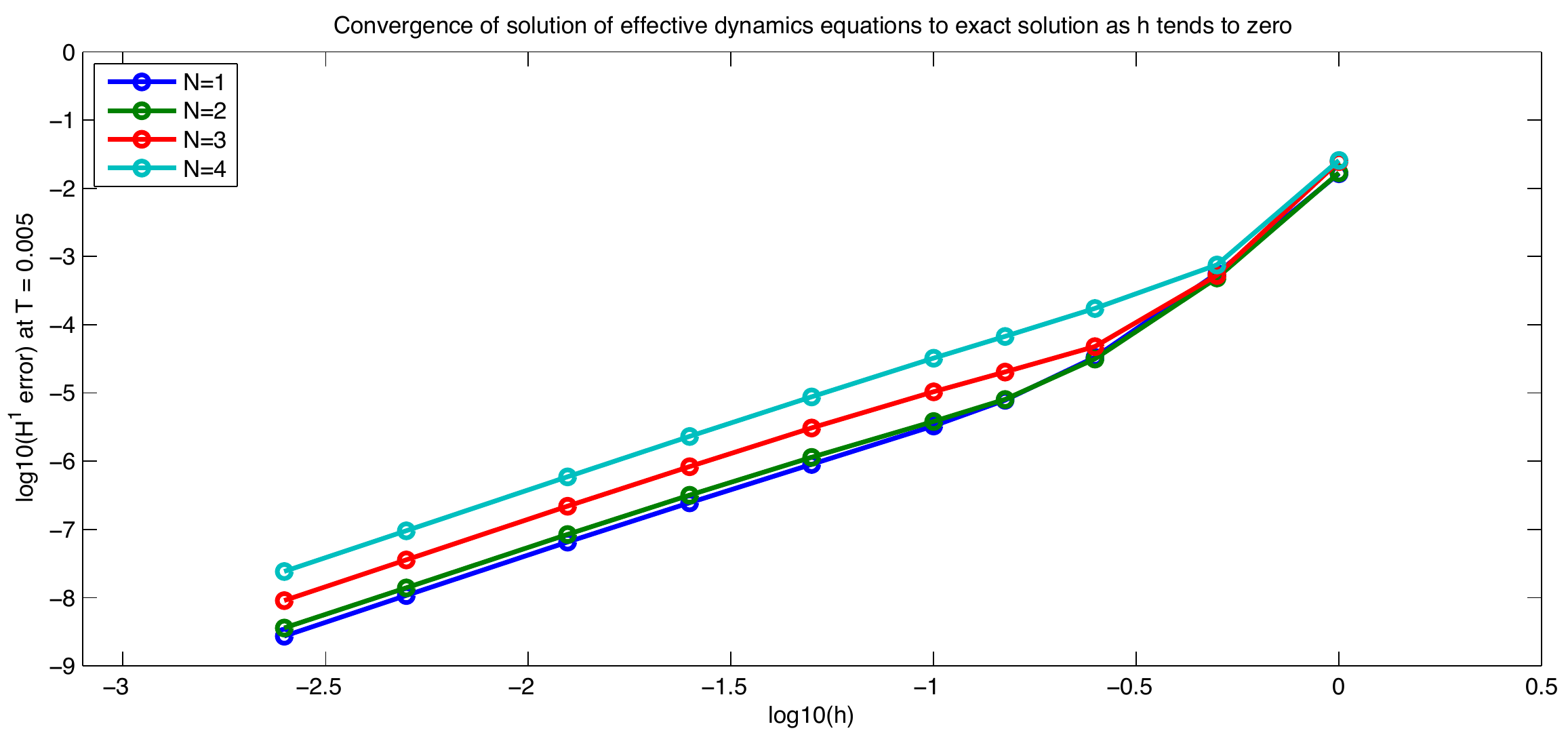}
\caption{A log-log plot of the $H^1$ error,
relative to the $H^1$ norm of the initial data,
 between the exact solution to the mKdV equation
(\ref{eq:mkdv}) and the $N$-soliton evolving according to the effective dynamics, 
as a function of $h$.
For smaller values of $h$, the slope of the lines approaches 2, 
in agreement with the theoretical upper bound on the error.
The theoretical upper bound has only been proven for $N = 1,2$,
but this figure gives evidence that it holds for all $N$. }
\label{fig:logherrmkdv}
\end{figure}

Similarly to \S \ref{subsec:htozero}, 
we have conducted a numerical study verifying that the $H^1$ error
is $\mathcal{O}(h^2)$ as  $h \to 0$ for multiple soliton initial conditions. 
See Figures \ref{fig:birdseyemkdv} and \ref{fig:logherrmkdv}.

\section{Numerical methods}
\label{sec:num}
We now describe the numerical methods we employ to compute
the Gross-Pitaevskii PDE (\ref{eq:GP}) 
and the ODE (\ref{eq:effdyn}) arising from the effective dynamics.
When comparing a solution of (\ref{eq:GP}) with a solution of (\ref{eq:effdyn}),
we refine our numerical solutions until the error between sucessive 
refinements of solutions to the same equation is several orders of magnitude
smaller than the error between solutions of the two equations.

Numerically solving the ODE arising from the effective dynamics (\ref{eq:effdyn})
 necessitates computing $q_N(x,a,v,\theta,\mu)$ and its derivatives with respect to
the parameters $a,v,\theta$,and $\mu$ efficiently. 
For this we note that an equivalent definition of $q_N$ in (\ref{eq:q_N})
\begin{equation}
\label{eq:altqndef}
q_N(x) = -\vec{1} M^{-1} \gamma \,,
\end{equation}
where $\vec{1}, M$, and $\gamma$ are as in (\ref{eq:q_N}).
Since $iM$ is Hermitian, $M^{-1}\gamma$ can be efficiently computed using 
the Cholesky factorization.
Differentiating $q_N$ numerically, for larger values of $N$,
is too costly.
Instead we used (\ref{eq:altqndef}) to obtain explicit formulas for the derivatives
of $q_N$ and again used Cholesky factorizations to efficiently compute them. 
With this we compute the integrand in 
(\ref{eq:V-N}), and then numerically integrate it using the trapezoidal method.
Once we can efficiently compute the RHS 
of the effective dynamics equations (\ref{eq:effdyn}),
the standard fourth-order Runge Kutta method was found to be suitable to solve to the ODE.

In order to solve the PDE (\ref{eq:GP}) we used a
Fourier spectral method to study the evolution on the 
numerical domain $(-\pi,\pi)$. 
This requires our solution $u(x,t)$ to be periodic in space, 
so we choose initial conditions such that $u(x,t)$ decays to zero,
to machine precision, before the endpoints $-\pi$ and $\pi$. 
Arbitrary initial data can be handled by either extending 
the numerical domain or rescaling the equation (see (\ref{eq:scalpha})).

One difficulty arises in that a non-trivial potential $W(hx)$ 
cannot be periodic for all $h$. 
However, the potential $W(hx)$ need not be periodic on $(-\pi, \pi)$ 
so long as the product $W(hx)u(x,t)$ is periodic on $(-\pi, \pi)$, 
which is achieved if $u(x,t)$ decays fast enough at the endpoints
$-\pi$ and $\pi$. 
If $W(x)$ is periodic on $(-\pi, \pi)$ 
and one wishes to consider a solution $u(x,t)$ that doesn't decay
before the endpoints $-\pi$ and $\pi$, 
the rescaling (\ref{eq:scalpha}) may be employed with $\alpha = h$:
\begin{equation}
\label{eq:scaleh}
\tilde x = hx\,, \ \ \tilde t = h^2 t\,, \ \ 
\tilde u(\tilde x,\tilde t)  \stackrel{\text{\tiny def}}{=} \frac{1}{h}u(x,t)\,. \ \
\end{equation}
Then if $\tilde u(\tilde x,\tilde t)$ satisfies (\ref{eq:GP})
with periodic potential $\tilde V(\tilde x) = W(\tilde x)/h^2$,
$ u(x,t)$ also satisfies (\ref{eq:GP}) with 
potential $V(x) = W(hx)$. 

This rescaling also makes it clear that as $h \to 0$, 
a soliton solution becomes sharper relative to the potential.
This requires higher resolution in order to apply our numerical method to solve
the PDE (\ref{eq:GP}), 
while the effective dynamics equations (\ref{eq:effdyn}) are unaffected. 
Indeed, our numerical experiments confirmed that increased computational effort was needed to
resolve the PDE as $h \to 0$, but not the effective dynamics ODE.

We now describe the method to solve a general solution $u(x,t)$ of (\ref{eq:GP}) 
on a periodic domain with periodic initial data and potential $V(x)$.
The Fourier modes $\hat{u}_k(t)$ of a solution $u(x,t)$ to (\ref{eq:GP}) 
evolve according to 
\begin{equation}
\label{eq:fourier}
\partial_t\hat{u}_k = -\frac{i}{2}k^2\hat{u}_k + i(\widehat{uv})_k,\quad v = |u|^2 - V
\end{equation} 
Discretizing space and replacing the Fourier Transform with the
Discrete Fourier Transform gives rise to a finite dimensional system of ODE, 
which we now represent in the general form 
\begin{equation}
\label{eq:LN}
u_t = \mathcal{L}u + \mathcal{N}(u)
\end{equation}
where $\mathcal{L}$ is a stiff linear transformation corresponding to the first term
of (\ref{eq:fourier}) (represented by a diagonal matrix in our case)
and $\mathcal{N}$ is a non-linear operator from the second term of (\ref{eq:fourier}).
To solve (\ref{eq:LN}) we compared the fourth order implicit-explicit (IMEX) method 
ARK4(3)6L[2]SA proposed by Kennedy and Carpenter \cite{KC}
with the exponential time differencing (ETD) method ETDRK4 proposed by 
Kassam and Trefethen \cite{KT}.
The IMEX scheme update formula is
\begin{equation}
u_{n+1} = u_n + \Delta t (b_1(k_1+l_1) + \cdots + b_s(k_s+l_s))\,,
\end{equation}
where $\Delta t$ is the time step and $k_i$ and $l_i$ are chosen such that
\begin{equation}
k_i = \mathcal{L}(u_n + \Delta t (A_{i1}k_1 + \cdots + A_{is}k_s + \hat A_{i1}l_1 + \cdots + \hat A_{is}l_s))
\end{equation}
\begin{equation}
l_i = \mathcal{N}(u_n + \Delta t (A_{i1}k_1 + \cdots + A_{is}k_s + \hat A_{i1}l_1 + \cdots + \hat A_{is}l_s))\,.
\end{equation}
Here $A, \hat A$ are $s\times s$ lower triangular matrices with 
$\hat A$ having zeros along its diagonal..
This allows us to solve for the $k_i$ and $l_i$ one stage at a time,
only inverting the diagonal linear operators 
$I - \Delta t A_{ii}\mathcal{L}$.
The implicit treatment of the $\mathcal{L}$ term mitigates 
the stiffness arising from the $k^2$ factor in
(\ref{eq:fourier}), while the $l_i$ can be computed
explicitly, so non-linear equations involving $\mathcal{N}$ need not be solved.

The ETD method, on the other hand, uses an exact formula for obtaining the next step 
$u_{n+1}$ from $u_n$ based on solving the linear portion exactly:
\begin{equation}
\label{eq:ETD}
u_{n+1} = e^{\mathcal{L}\Delta t}u_n + e^{\mathcal{L}\Delta t}
\int_0^{\Delta t} e^{-\mathcal{L}\tau} \mathcal{N}(u(t_n + \tau),t_n + \tau)\, d\tau
\end{equation}
The integral in (\ref{eq:ETD}) can then be numerically approximated using matrix exponents of
$\mathcal{L}$ and evaluations of $N$. 
Thus as with the IMEX method, we do not need to solve non-linear equations, and computations
involving $\mathcal{L}$ (namely computing $e^{\Delta t \mathcal{L}}$) are efficient because 
$\mathcal{L}$ is diagonal. Stiffness is mitigated by solving the
linear portion of (\ref{eq:LN}) exactly.
We found that the ETDRK4 scheme computed a solution of a desired accuracy nearly twice
as fast as the ARK4(3)6L[2]SA scheme. 
While the ARK4(3)6L[2]SA scheme had a slightly smaller error rate per step,
more computations per step made it significantly less efficient. 
Neither method demonstrated any instability in the range of step sizes required 
for our solutions.
Below, we plot the convergence of the two schemes as the timestep goes to zero. 
To obtain the results in the figures, we used the same potential 
and initial data as in \S \ref{subsec:ancs}: 
$W(x) = -100\, \text{sech}^2 (5x) + 10x^4\,$.

\begin{figure}[t]
\includegraphics[scale=.5]{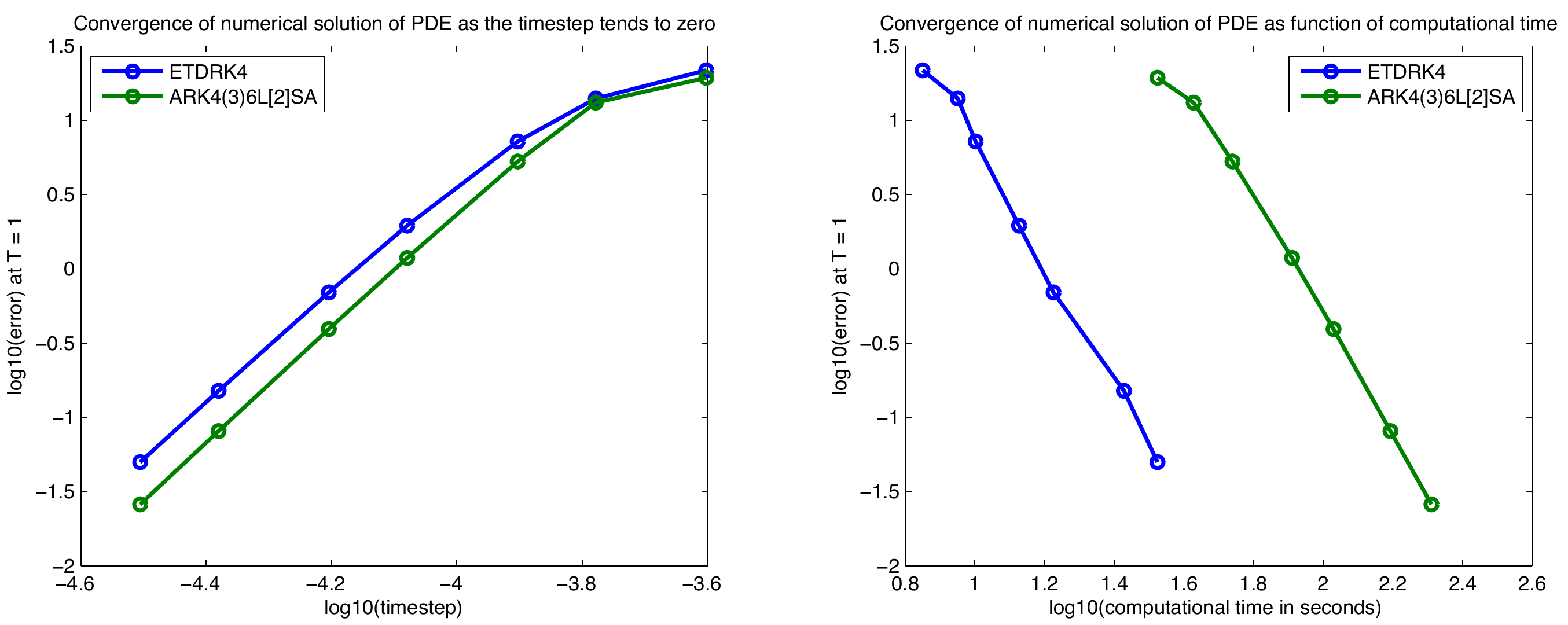}
\caption{Log-log plots of the convergence of the fourth order schemes
ETDRK4 and ARK4(3)6L[2]SA as a function of the timesteps and computational time, respectively.
ARK4(3)6L[2]SA is slightly more efficient per timestep, but ETDRK4 is significantly more
computationally efficient.}
\label{fig:logetdimex}
\end{figure}

\vspace{2 mm}




\subsection*{Acknowledgements}
The author was supported through the National Science Foundation through 
grant DMS-0955078 and by the Director, 
Office of Science, Computational and Technology Research, 
US Department of Energy under Contract DE-AC02-05CH11231.
The author would like to thank Jon Wilkening and Maciej Zworski 
for helpful discussions and comments.
\bibliographystyle{amsplain}
\bibliography{refs2}

\end{document}